\documentclass[]{article}

\usepackage[]{babel}
\usepackage[T1]{fontenc}
\usepackage[utf8]{inputenc}
\usepackage{amsmath}
\usepackage{amssymb}
\usepackage{amsthm}
\usepackage{graphicx}
\usepackage{enumerate}
\usepackage{cite}

\makeatletter

\@addtoreset{equation}{section}
\makeatother

\begin{document}
	
	\newcommand{\Gam}{\Gamma(1)}
	\newcommand{\C}{\mathbb{C}}
	\newcommand{\R}{\mathbb{R}}
	\newcommand{\Q}{\mathbb{Q}}
	\newcommand{\Qcal}{\mathcal{Q}}
	\newcommand{\Scal}{\mathcal{S}}
	\newcommand{\Z}{\mathbb{Z}}
	\newcommand{\N}{\mathbb{N}}
	\newcommand{\F}{\mathbb{F}}
	\newcommand{\M}{\mathcal{M}}
	\newcommand{\s}{\mathcal{S}}
	\newcommand{\h}{\mathfrak{h}}
	\newcommand{\Lcal}{\mathcal{L}}
	\newcommand{\achi}{\mathfrak{a}_\chi}
	\newcommand{\abcd}{\begin{pmatrix}
			a&b\\c&d
	\end{pmatrix}}
	\newcommand{\matp}{\gamma_{q^A}}

	\newcommand{\ddx}[1]{\frac{\partial #1}{\partial x}}
	\newcommand{\haar}[1]{\frac{d#1}{#1}}
	\newcommand{\ol}[1]{\overline{#1}}
	\newcommand{\wt}[1]{\widetilde{#1}}
	\newcommand{\G}[1]{\Gamma(#1)}
	\newcommand{\g}[1]{\Gamma_0(#1)}
	\newcommand{\gup}[1]{\Gamma_0^0(#1)}
	\newcommand{\GL}[1]{\mathrm{GL}(#1,\Z)}
	\newcommand{\SL}[1]{\mathrm{SL}(#1,\Z)}
	\renewcommand{\Im}[1]{\textup{Im}(#1)}
	\renewcommand{\Re}[1]{\textup{Re}(#1)}
	
	\newcommand{\footremember}[2]{%
		\footnote{#2}
		\newcounter{#1}
		\setcounter{#1}{\value{footnote}}%
	}
	\newcommand{\footrecall}[1]{%
		\footnotemark[\value{#1}]%
	} 
	
	\newtheorem{thm}{Theorem}[section]
	\newtheorem{definition}[thm]{Definition}
	\newtheorem{conj}[thm]{Conjecture}
	\newtheorem{lemme}[thm]{Lemma}
	\newtheorem{ex}[thm]{Example}
	\newtheorem{cor}[thm]{Corollary}
	\newtheorem{exo}[thm]{Exercise}
	\newtheorem{prop}[thm]{Proposition}
	\newtheorem{rem}[thm]{Remark}

\title{Non vanishing of products of twisted $\mathrm{GL}{(3)}$ L-functions}
\author{Robin Frot
\\
	\small{Alfr\'ed R\'enyi Institute of Mathematics, POB 127, Budapest H-1364,	Hungary}\\
\small{MTA R\'enyi Int\'ezet Lend\"ulet Automorphic Research Group}}

\maketitle

\tableofcontents


\begin{abstract}
	In this paper, we prove that if the Fourier coefficients of a $\SL{3}$ Hecke--Maa\ss\ cusp form $\pi$ are not too correlated with additive characters, then there exists infinitely many Dirichlet characters such that
	\begin{align*}
		L\left(\frac{1}{2},\pi\otimes\chi\right)L\left(\frac{1}{2},\chi\right)\neq 0.
	\end{align*}
	To prove this result, we compute the first twisted moment of these $L$ function averaged over a well chosen set of conductors.
\end{abstract}

\section{Introduction}

The main objective of this paper is to prove a simultaneous non-vanishing result of products twisted Hecke--Maa\ss\ L-functions and Dirichlet L-functions. We manage to do so by computing the average of the first moment over the moduli. Note that similar results for $\GL{2}$ Hecke--Maa\ss\ forms has been done by Das and Khan in \cite{das_simultaneous_2014}. Indeed, if $f$ is an even Maa\ss\ form and $q$ a prime number, one have
\begin{align*}
	\frac{1}{\phi(q)}\sideset{}{^*}\sum_{\chi [q]}L\left(\frac{1}{2},f\otimes\chi\right)\ol{L\left(\frac{1}{2},\chi\right)} = L(1,f)+ 0\left(q^{-\frac{1}{64}+\epsilon}\right)
\end{align*}
where the star means that the sum is taken over even primitive characters.

In the case of $\GL{3}$ L-functions, a similar work has been done by Munshi and Sengupta without the Dirichlet $L$-function in \cite{munshi_determination_2014} if we do not take the product with the associated Dirichlet $L$-function. However, an additional average over the conductors $q$ is needed to apply a method developed by Luo in \cite{luo2005nonvanishing}. They proved that if $\pi$ is a Hecke--Maa\ss\ cusp form for $\SL{3}$, then 
\begin{align*}
	\frac{1}{|\mathcal{Q}|}\sum_{q\in\mathcal{Q}}\frac{1}{\phi(q)}\sum_{\chi[q]}L\left(\frac{1}{2},\pi\otimes\chi\right) = 1 + O(|\Qcal|^{-\eta+\epsilon})
\end{align*}
for some $\eta>0$, any $\epsilon>0$ and well chosen set of integers $\Qcal$.
We aim at estimating the sums
\begin{align}
	\Scal_\Qcal :=\frac{1}{|\mathcal{Q}|}\sum_{q\in\mathcal{Q}}\frac{1}{\phi(q)}\sum_{\chi[q]}L\left(\frac{1}{2},\pi\otimes\chi\right)L\left(\frac{1}{2},\ol \chi\right).
\end{align}
Using the methods of Luo in \cite{luo2005nonvanishing} in the $GL(4)$ case, one can get
\begin{align*}
	\frac{1}{|\mathcal{Q}|}\sum_{q\in\mathcal{Q}}\frac{1}{\phi(q)}\sum_{\chi[q]}L\left(\frac{1}{2}+\alpha,\pi\otimes\chi\right)L\left(\frac{1}{2}+\alpha, \chi\right) \gg 1+O\left(|\Qcal|^{-2\Re \alpha + \epsilon}\right)
\end{align*}
when $\alpha$ is a complex number of positive real part. The estimate for the moment fails exactly at the critical line and we need additional cancellation to deal with the non-vanishing at the central value.
We prove here they following result.
\begin{thm}\label{mainThm}
	Let $\pi$ be an Hecke--Maa\ss\ cusp form for $\SL{3}$ and let $A(n,1)$ be its Fourier--Whittaker coefficients. Assume that the satisfy the following:
	\begin{align}\label{cond}
		\sum_{n\leq X}A(n,1)e(nx)\ll X^{\frac{1}{2}+a+\epsilon}
	\end{align}
uniformly in $x$, for some $a<\frac{1}{6}$ and any $\epsilon>0$. Then, for a good choice of set $\Qcal$, there exists $\eta = \eta(\delta)>0$ such that
\begin{align}
	\frac{1}{|\mathcal{Q}|}\sum_{q\in\mathcal{Q}}\frac{1}{\phi(q)}\sum_{\chi[q]}L\left(\frac{1}{2},\pi\otimes\chi\right)L\left(\frac{1}{2},\ol \chi\right) = L(1,\pi)+O\left(|\Qcal|^{-\eta + \epsilon}\right)
\end{align}
for any $\epsilon>0$.
\end{thm}
One can note that an additional condition on the cusp form is needed to improve Luo's result. Indeed, every power saving in Luo's method is only coming form a careful analysis of characters sums. 

The additional condition is given here in the form of a non correlation between the Fourier coefficients and additive characters. In \cite{miller2006automorphic}, Miller proved that one can always choose $a = 1/4$ in \eqref{cond}. However, it is conjectured that the optimal constant $a = 0$ can be reached. As proved in Miller's article this might either be done by computing directly the sums 
\begin{align*}
	\sum_{n\leq T} A(n,1)e(\alpha n)
\end{align*}
or by estimating the periods of the Maa\ss\ form. In any case, this conjecture is consider very hard. Hopefully, we only need a slightly weaker version in which we do not have to assume the full square root cancellation.
On can note that if the condition \eqref{cond} on the $\SL{3}$ Maa\ss\ forms always holds, this result complete our knowledge of (averaged) first twisted moment at $1/2$ of product of $L$ function of total degree $4$. The case of four Dirichlet $L$ functions is treated by Young in \cite{young}, and the case of two $\GL{2}$ L-functions is treated in \cite{blomer2017moments}.

\paragraph{Acknowledgements:}

This work is the conclusion of the author's PhD thesis funded by the \'ENS Lyon and the University of Lille. The authors want to thank his supervisors Gautami Bhowmik and Nicole Raulf, without whom, this work would not have been possible. This article was written during a post-doc at the R\'enyi Institue in Budapest.

\section{Background on Maa\ss\ forms and main tools}
 
 \subsection{Dirichlet characters}
 
 As usual, we define the Dirichlet characters modulo $n$ as morphism
 \begin{align*}
 	\chi: \left(\Z/n\Z\right)^*\to \C^*
 \end{align*}
and we extend them $n$-periodically to $\Z$ with $\chi(k) = 0$ if $(n,k)\neq 1$. If $\chi$ is a Dirichlet character modulo $n$, we define its conductor $q$ as the smallest integer $q|n$ such that we can factor $\chi$ through the projection map
\begin{align*}
	\left(\Z/n\Z\right)^*\to \left(\Z/q\Z\right)^*.
\end{align*}
If $\chi$ is a Dirichlet character modulo $q$ with conductor $q$, we say that $\chi$ is primitive.
A Dirichlet character is said to be even (resp. odd) if  $\chi(-1)=1$ (resp. $\chi(-1)=-1$).
 For any primitive Dirichlet character modulo $q$, we can define the normalized Gauss sum 
\begin{align}
	\epsilon(\chi) := \frac{1}{\sqrt{q}}\sum_{a=0}^{q-1}e_q(a)\chi(a)
\end{align}
where
\begin{align*}
	e_q(n) := e\left(\frac{n}{q}\right) = \exp\left(\frac{2i\pi n}{q}\right).
\end{align*}
We have the following orthogonality relation for even Dirichlet characters.
\begin{prop}\label{characterSums}
	Let $q_1\neq q_2$ be two prime numbers and $k$ be a positive integer and $n$ such that $(n,q)=1$. Then we have 
	\begin{align}
		\frac{2}{\phi(q)}\sideset{}{^*}\sum_{\substack{\chi_1 [q_1]\\\chi_2[q_2]}}\chi_1\chi_2(n)& = \mathbf{1}_{n=\pm1[q_1q_2]}-\frac{\mathbf{1}_{n=\pm1[q_2]}}{\phi(q_1)}-\frac{\mathbf{1}_{n=\pm1[q_1]}}{\phi(q_2)}+\frac{2}{\phi(q_1q_2)}\\
		\frac{2}{\phi(q)}\sideset{}{^*}\sum_{\substack{\chi_1 [q_1]\\\chi_2[q_2]}}\epsilon(\ol {\chi_1\chi_2})\chi_1\chi_2(n) &= \frac{1}{\sqrt{q}}\sum_{\alpha\in\{\pm1\}}\sum_{d|q}\frac{1}{\phi(d)}e_{\frac{q}{d}}(\alpha mn)
	\end{align}
where the stars over the sums are indicating they are taken over primitive characters such that $\chi_1\chi_2$ is an even character of conductor $q_1q_2$.
\end{prop}

To complete this section, we extract a lemma from \cite{luo2005nonvanishing}. This lemma is Luo's key idea for his computation of the first twisted moment of $\mathrm{GL}(n)$ L-functions.
Let us define the set 
\begin{align}
	\Qcal := \left\{q=q_1q_2, \ q_1,q_2\text{ primes}, \ Q_1<q_1<2Q_1,\ Q_2<q_2<2Q_2\right\}.
\end{align}
Moreover, let us write $Q = Q_1Q_2$. We will also assume that $c<Q_1<Q_2$ for a large enough constant $c$. Note that if $Q$ is large enough,
\begin{align}
	|\mathcal{Q}|\gg Q^{1-\epsilon}.
\end{align}
Let us define
\begin{align}\label{H}
	H(y) := \frac{1}{\pi}(1+y^{2})^{-1}.
\end{align}
and for any integer $k\geq1$,
\begin{align}\label{Bk}
	B_k(s,m):=\frac{1}{|\mathcal{Q}|}\sum_{q\in\Qcal}\frac{q^{s}}{\phi(q)}\sideset{}{^*}\sum_{\substack{\chi [q_1]}}\epsilon(\ol{\chi})^k\chi(m).
\end{align}
The star above the sum signifies that we sum over the characters such that $\chi_1\chi_2$ is even and primitive.
We have the following
\begin{lemme}[Luo]\label{Luo}
	For any real number $y\gg 1$ and integer $k\geq 1$, we have
	\begin{align*}
		\sum_{m\in\Z}H\left(\frac{m}{y}\right)|B_k(s,m)|^2\ll Q^{\Re s}\left(y Q^{-2+\epsilon}Q_1^2 + Q^\epsilon Q_2^{-\frac{1}{2}}+Q^{-1+\epsilon}Q_2\right)
	\end{align*}
	Moreover, by taking $Q_1\asymp Q^{\delta}$ with $\delta<1/2$, we get
	\begin{align}
		\sum_{m\in\Z}H\left(\frac{m}{y}\right)|B_k(s,m)|^2\ll Q^{\Re s}\left(y Q^{-2+2\delta+\epsilon} + Q^{-\frac{1}{2}+\frac{\delta}{2}+\epsilon}+Q^{-\delta+\epsilon}\right).
	\end{align}
\end{lemme}

This Lemma can be extracted from the series of calculation made by Luo in \cite{luo2005nonvanishing} (after formula (4)). In our case, the central character $\omega_\pi$ of $\pi$ is trivial.

 \subsection{Hecke--Maa\ss\  forms for $\SL{3}$}
 
 The reader could refer to \cite{goldfeld2006automorphic} for a more in depth exposition of the theory of Hecke--Maa\ss\ forms for $\SL{N}$. Except if it is mentioned explicitly, this section will refer to Goldfeld's book.
 
 Let us give some background about the theory of Hecke--Maa\ss\ forms. Define the generalized Poincaré half-plane as
 \begin{align}
 	\h^3:= \mathrm{GL}(3,\R)/\left<O(3,\R),\R^\times\right>
 \end{align}
where $\R^\times$ is assimilated to the center of $\mathrm{GL}(3,\R)$.
 The group $\SL{3}$ acts discretely on $\h^3$ by left multiplication. Finally let us write $\mathcal{D}^3$ to be the center of the universal enveloping algebra of $\mathfrak{gl}(3,\R)$.
 
 Every element $z\in\h^3$ can be written in the form
 \begin{align*}
 	z=x\cdot y = \begin{pmatrix}
 		1&x_1&x_3\\
 		0&1&x_2\\
 		0&0&1
 	\end{pmatrix}\cdot\begin{pmatrix}
 	y_1y_2&0&0\\
 	0&y_1&0\\
 	0&0&1
 \end{pmatrix}
 \end{align*}
where $x_1,x_2,x_3\in\R$ and $y_1,y_2>0$. We can define the Haar measure on $\h^3$ by $$dz:=\frac{dx_1dx_2dx_3dy_1dy_2} {y_1^3y_2^3}$$ which is invariant under left multiplication by $\mathrm{SL}(3,\R)$. The invariance of the measure allows us to define it on the quotient $\SL{3}\backslash\h^3$. Note that this quotient is of finite volume. We can also define the Peterson hermitian product by 
\begin{align}
	\left<f,g\right> = \int_{\SL{3}\backslash\h^3}f(z)\ol{g(z)}dz.
\end{align}

Finally let us write $\mathcal{D}^3$ to be the center of the universal enveloping algebra of $\mathfrak{gl}(3,\R)$.
 \begin{definition}
 	A Maa\ss\ cusp for for $\SL{3}$ is a function in $\mathcal{L}^2\left(\SL{3}\backslash \h^3\right)$
 	such that
 	\begin{itemize}
 		\item $f$ is an eigenfunction of all differential operators in $\mathcal{D}^3$,
 		\item for every unitary subgroup
 		\begin{align*}
 			U_1 = \left\{\begin{pmatrix}
 				1&*&*\\0&1&*\\0&0&1
 			\end{pmatrix}\right\},\
 		U_2 = \left\{\begin{pmatrix}
 			1&0&*\\0&1&*\\0&0&1
 		\end{pmatrix}\right\},\
 	U_3 = \left\{\begin{pmatrix}
 		1&*&*\\0&1&0\\0&0&1
 	\end{pmatrix}\right\},
 		\end{align*}
 	we have 
 	\begin{align*}
 		\int_{\SL{3}\cap U_i\backslash U_i}f = 0.
 	\end{align*}
 	\end{itemize} 
 \end{definition}

The Hecke operators for $\SL{3}$ are defined in the following way. Let $n\geq1$ be an integer, the n-th Hecke operator $T_n$ is defined by
\begin{align*}
	T_n f(z) = \sum_{\substack{abc = n\\ 0\leq b_1< b\\ 0\leq c_1,c_2< c}}f\left(\begin{pmatrix}
		a& b_1 & c_1\\ 0 &b&c_2 \\ 0&0&c
	\end{pmatrix}\cdot z\right).
\end{align*}
Note that they commute with every operators in $\mathcal{D}^3$
 An Hecke--Maa\ss cusp for $\SL{3}$ is a Maa\ss\ cusp for that is also an eigenfunction of all the Hecke operators.
 
 Every Hecke--Maa\ss\ cusp form $\pi$ has a Fourier--Whittaker expansion of the following shape:
 \begin{align}
 	\pi(z) = \sum_{\gamma\in U_2(\Z)\backslash\SL{2}}&\sum_{m_1>0}\sum_{m_2\neq0}\frac{A(m_1,m_2)}{|m_1m_2|}\\&\times W_\nu\left(\begin{pmatrix}
 		|m_1m_2|&0&0\\0&m_1&0\\0&0&1
 	\end{pmatrix}\cdot\begin{pmatrix}
 	\gamma&0\\0&1
 \end{pmatrix}\cdot z,\phi_{(1,\frac{m_2}{|m_2|})}\right)\nonumber
 \end{align}
where $W_\nu$ is the Jacquet--Whittaker function of type $\nu$.
Assume that $\pi$ is normalized such that $A(1,1) = 1$. Then its coefficients must satisfy the following
\begin{align}
	\ol{A(m,n)} &= A(n,m),\\
	A(m,-n)&= A(m,n),\\
	A(m_1m_2, n_1n_2) &= A(m_1,n_1)A(m_2,n_2) \text{ if }(m_1n_1,m_2n_2)=1\label{Hecke}\\
	A(m,n)&= \sum_{d|(m,n)}\mu(d)A\left(\frac{m}{d},1\right)A\left(1,\frac{n}{d}\right).
\end{align}

From these relation it is possible to define the standard $L$-function associated to a Hecke--Maa\ss\ cusp form for $\SL{3}$. If $\chi$ a Dirichlet we also give the definition of the twisted $L$-function. For $\Re s>1$, we define
\begin{align}
	L(s,\pi)&:= \sum_{n\geq1}\frac{A(n,1)}{n^s} = \prod_{p}\left(1 - \frac{A(p,1)}{p^s}+ \frac{A(1,p)}{p^{2s}}-\frac{1}{p^{3s}}\right)^{-1}\\
	L(s,\pi\otimes \chi)&:= \sum_{n\geq1}\frac{A(n,1)\chi(n)}{n^s}\\&= \prod_{p}\left(1 - \frac{A(p,1)\chi(p)}{p^s}+ \frac{A(1,p)\chi(p^2)}{p^{2s}}-\frac{\chi(p^3)}{p^{3s}}\right)^{-1}
\end{align}
 
 These $L$ functions have an analytic continuation to $\C$ with no poles and they satisfy a functionnal equation. Assume that the character $\chi$ is even and primitive of conductor $q$, then 
 \begin{align}
 	\Lambda(s,\pi\otimes\chi) := q^\frac{3s}{2}G_\pi(s)L(s,\pi\otimes\chi) = \epsilon(\chi)^3\Lambda(1-s,\wt\pi\otimes\ol\chi)
 \end{align}
where
\begin{align}
	G_\pi(s) =\pi^{-\frac{3s}{2}} \prod_{i=1}^3\Gamma\left(\frac{s+\mu_i}{2}\right)
\end{align}
and $\wt \pi$ is the dual Maa\ss\ form of $\pi$.
The coefficients $\mu_i$ are the Langlands parameters of $\pi$.

Thanks to the work of Kim and Sarnak in \cite{kim_functoriality_2003}, we know that for all $i$,
\begin{align}
	|\Re{\mu_i}| < 5/14.
\end{align}

 \subsection{Approximate functional equations}
 
 Since our goal is to study the above $L$ functions outside of the domain of convergence, we will need to express these functions in term of the coefficients. In order to do so, we apply the approximate functional equations. They are proved for instance in \cite{iwaniec2004analytic}, Chapter 5.3.
 
 \begin{lemme}[Approximate functional equations]
 	Let $\pi$ be a Hecke--Maa\ss\ cusp form for $\SL{3}$ and $\chi$ be an even primitive Dirichlet character. We have 
 	\begin{align}
 		L\left(\frac{1}{2},\pi\otimes\chi\right)L\left(\frac{1}{2},\ol\chi\right) =& \sum_{m\geq1}\frac{A(m,1)\chi(m)\ol{\chi}(n)}{\sqrt {mn}}V\left(\frac{mn}{q^{2-\eta}}\right)\\& + \epsilon(\chi)^2\sum_{m\geq1}\frac{A(1,m)\chi(n)\ol\chi(m)}{\sqrt {mn}}V\left(\frac{mn}{q^{2+\eta}}\right),\label{apeq1}\\
 	\end{align}
 where
 \begin{align*}
 	V(y)=\frac{1}{2\pi i}\int_{(2)}y^{-s}\frac{G\left(\frac{1}{2}+s\right)}{G\left(\frac{1}{2}\right)}\haar{s}\\
 \end{align*}
with $G(s) = G_\pi(s)\Gamma\left(\frac{s}{2}\right)$.
Moreover, the function $V$ satisfies
\begin{align}\label{boundV}
	V(y) = \left\{\begin{aligned}
		&1+O(y^{\frac{1}{7}+\epsilon}),\\
		&O_A(y^{-A}).
	\end{aligned}\right.
\end{align}
 \end{lemme}
 
 One can note that the exponent $1/7 = 1/2 - 5/14$ in \eqref{boundV} can be replaced by $1/2$ if we assume that the form is tempered (i.e. $\Re{\mu_i}=0$).
 
 \subsection{Bounds for the Fourier coefficients}
 
 For $\SL{3}$ Hecke--Maa\ss\ cusp forms, one has the following bounds for the Fourier coefficients. 
 \begin{align}
 	&|A(n,1)|< n^{\theta+\epsilon}\label{RP},\\
 	&\sum_{n\leq X} |A(n,1)|^2\ll X\label{RS},\\
 	&\sum_{n\leq X} A(n,1)\ll X^{\frac{3}{4}+\epsilon}\label{Miller}.
 \end{align}
 
 All these bounds hold for $\epsilon>0$.
 In the first bound, $\theta$ can be chosen to be as small as $5/14$ due to Kim and Sarnak \cite{kim_functoriality_2003}. The second bound is due to the study of the Rankin--Selberg $L$ function and can be found in \cite{goldfeld2006automorphic}. Finally, the last bound is due to Miller in \cite{miller_coef}. Miller's bound will not be enough for our purpose and we will need to assume that
 \begin{align*}
 	\sum_{n\leq X} A(n,1)\ll X^{\frac{1}{2}+a+\epsilon}
 \end{align*}
with $a<1/6$. A folkloric conjecture states that we can chose $a$ to be $0$.

\section{Proof of the main theorem}

Let us fix $Q\geq 1$ a sufficiently large integer and $0<\delta<1/2$ 
For $i=1,2$, let us define the sets
\begin{align*}
	\Qcal_i = \{q_i\in [Q_i,2Q_i], \ q_i \ is\ prime\}
\end{align*}
where $Q_1 = Q^\delta$ and $Q_2 = Q^{1-\delta}$.
 Also define $\Qcal\subset[Q,4Q]$ to be set following set of integer
\begin{align}
	\Qcal =\left\{q = q_1q_2, q_1\in\Qcal_1,q_2\in\Qcal_2, \right\}
\end{align}

Let $\pi$ be a Hecke--Maa\ss\ cusp form for $\SL{3}$ and define the averaged first moment as
\begin{align}
	\Scal_\Qcal := \frac{2}{|\Qcal|}\sum_{q\in \Qcal}\frac{1}{\phi(q)}\sideset{}{^*}\sum_{\chi [q]}L\left(\frac{1}{2},\pi\otimes\chi\right)L\left(\frac{1}{2},\ol\chi\right) \label{moment}
\end{align}
 where the star indicates that we sum over the even primitive Dirichlet characters modulo $q$.
 
In this section, we prove the following theorem
\begin{thm}\label{mainThm2}
	Let $\pi$ be Hecke--Maa\ss\ cusp form for $\SL{3}$ such that its Fourier coefficients satisfy the following
	\begin{align}
		\sum_{n\leq X} A(n,1)e(nx)\ll X^{\frac{1}{2}+a+\epsilon}\label{hyp}
	\end{align}
 	uniformly in $x$ for some $a<1/6$. With the above notation, we have for any $\epsilon>0$, $\delta<1/4-3a/2$,
 	\begin{align}
		\Scal_\Qcal = L(1,\pi) + O\left(Q^{\frac{3}{4}a-\frac{1}{8}+\frac{\delta}{2}+\epsilon}+Q^{-\frac{\delta}{2}+\epsilon}\right)
 	\end{align}
 where the constant only depends on $\epsilon$ an $\pi$.
\end{thm}
One can note that a similar estimates will hold for $a\geq1/6$. However, in this case, the "error" term may be arbitrarily large and cancel out with the main term.
The Theorem \ref{mainThm2} implies the following corollary.
\begin{cor}\label{cor}
	Let $\pi$ be Hecke--Maa\ss\ cusp form for $\SL{3}$ satisfying \eqref{hyp}. There exists infinitely many Dirichlet characters such that
	\begin{align*}
		L\left(\frac{1}{2},\pi\otimes\chi\right)L\left(\frac{1}{2},\chi\right)\neq 0.
	\end{align*}
\end{cor}
Finally, one can notice that the proof will not rely on the fact that we are on the central point, we can replace $1/2$ by $1/2+it$ and have the same conclusion.

\subsection{Applying the approximate functional equations}
In order to access the coefficients of the $L$ functions involved in $\Scal_\Qcal$ at $1/2$, we first need to apply the approximate functional equations. By writing the product of \eqref{apeq1} to be determined later, we get
\begin{align*}
	\Scal_\Qcal = S_0+S_1
\end{align*}
where
\begin{align}
	&\label{S_0} S_0 := \frac{2}{|\Qcal|}\sum_{q\in \Qcal}\frac{1}{\phi(q)}\sideset{}{^*}\sum_{\chi [q]}\sum_{m,n\geq1}\frac{A(m,1)\chi(m)\ol\chi(n)}{\sqrt{mn}}V\left(\frac{mn}{q^{2-\eta}}\right),\\	&\label{S_1} S_1 := \frac{2}{|\Qcal|}\sum_{q\in \Qcal}\frac{1}{\phi(q)}\sideset{}{^*}\sum_{\chi [q]}\epsilon(\chi)^2\sum_{m,n\geq1}\frac{A(1,m)\ol\chi(m)\chi(n)}{\sqrt{mn}}V\left(\frac{mn}{q^{2+\eta}}\right).
\end{align}
The objective is to compute the main term from $S_0$, which is usual in this kind of problem since the character sum give a congruence relation between $m$ and $n$. Finally we will use the additional condition \eqref{hyp} on the coefficient of $\pi$ and Luo's Lemma to bound the sums $S_1$.

\subsection{Computation of the main term in $S_0$}

After computing the sum over the characters using the first part of Proposition \ref{characterSums}, we are left with 
\begin{align*}
	S_0 = \frac{1}{|\Qcal|}\sum_{(q_1,q_2)\in \Qcal_1\times\Qcal_2}& \sum_{\substack{m,n\geq1\\(mn,q_1q_2)=1 }}\frac{A(m,1)}{\sqrt{mn}}V\left(\frac{mn}{q^{2-\eta}}\right)\\
	&\times \left(\mathbf{1}_{m=\pm n[q_1q_2]}-\frac{\mathbf{1}_{m=\pm n[q_2]}}{\phi(q_1)}-\frac{\mathbf{1}_{m=\pm n[q_1]}}{\phi(q_2)}+\frac{2}{\phi(q_1q_2)}\right)
\end{align*}
First let us compute the sum
\begin{align*}
	\frac{1}{|\Qcal|}\sum_{q_1q_2\in \Qcal}& \sum_{\substack{m\equiv \pm n [q_1q_2]\\(mn,q_1q_2)=1 }}\frac{A(m,1)}{\sqrt{mn}}V\left(\frac{mn}{q^{2-\eta}}\right).
\end{align*}
As usual, we will get the main term from the diagonal terms $m=n$. They contribute to 
\begin{align*}
	\frac{1}{|\Qcal|}&\sum_{q_1q_2\in \Qcal}\sum_{(m,q_1q_2)=1}\frac{A(m,1)}{m}V\left(\frac{m^2}{q^{2- \eta}}\right)\\
	=& \frac{1}{(2\pi i)|\Qcal|}\sum_{q_1q_2\in \Qcal}\int_{(2)}q^{(2-\eta)s}L_{q_1,q_2}(1+2s)\frac{G\left(\frac{1}{2}+s_1\right)}{G\left(\frac{1}{2}\right)}\haar{s}
\end{align*}
where we used the definitions of the function $V$ and where the function $L_{q_1,q_2}$ is given by
\begin{align*}
	L_{q_1,q_2}(s):= \prod_{i=1}^{2}\left(1-\frac{A(q_i,1)}{q_i^s}+\frac{A(1,q_i)}{q_i^{2s}}-\frac{1}{q_i^{3s}}\right)L(s,\pi).
\end{align*} 
By shifting the line of integration to the left of $0$, we get the pole at $s=0$. Moreover we can shift the line up to $\Re{s} = -1/7+\epsilon = -1/2 + 5/14+\epsilon
$ to avoid meeting a pole of $G$. This implies that the contribution of the diagonal terms is 
\begin{align*}
	\frac{1}{|\Qcal|}&\sum_{q_1q_2\in \Qcal}L_{q_1,q_2}(1) + O(Q^{\frac{\eta-2}{7}+\epsilon}).
\end{align*}
Using the bound \eqref{RP} with $\theta = 5/14$, we have \begin{align*}
	\left(1-\frac{A(q_i,1)}{q_i}+\frac{A(1,q_i)}{q_i^2}-\frac{1}{q_i^{3}}\right) = 1 + O(Q_i^{-\frac{9}{14}+\epsilon}).
\end{align*}
Finally, the diagonal terms contribute to \begin{align}\label{Diag}
	L(1,\pi)+O(Q^{\frac{\eta-2}{7}+\epsilon}+Q^{-\frac{9\delta}{14}+\epsilon}).
\end{align}
We can now deal with the off-diagonal terms. Using the bound \eqref{boundV}, we can cut the tail of the sum for a negligible cost
\begin{align*}
	\frac{1}{|\Qcal|}\sum_{q_1q_2\in \Qcal}& \sum_{\substack{m\equiv \pm n [q_1q_2]\\(mn,q_1q_2)=1\\n\neq m }}\frac{A(m,1)}{\sqrt{mn}}V\left(\frac{mn}{q^{2-\eta}}\right)\\&\ll  \frac{1}{|\Qcal|}\sum_{q_1q_2\in \Qcal} \sum_{\substack{m\equiv \pm n [q_1q_2]\\mn<Q^{2-\eta+\epsilon}\\n\neq m }}\frac{|A(m,1)|}{\sqrt{mn}} + O(Q^{-100}).\\
	&\ll \frac{1}{|Q|}\sum_{\substack{m\neq n\\mn<Q^{2-\eta+\epsilon}}}\frac{|A(m,1)|}{\sqrt{mn}}\sum_{\substack{q\in\Qcal\\ q|(m\pm n)}}1+O(Q^{-100})\\
	&\ll Q^{-1+\epsilon}\sum_{m<Q^{2-\eta+\epsilon}}|A(m,1)|m^{-\frac{1}{2}+\epsilon}\sum_{n<\frac{Q^{2-\eta+\epsilon}}{m}}n^{-\frac{1}{2}+\epsilon}+O(Q^{-100})\\
	&\ll Q^{-\frac{\eta}{2}+\epsilon}\sum_{m<Q^{2-\eta+\epsilon}}|A(m,1)|m^{-1+\epsilon}+O(Q^{-100})\\
	&\ll Q^{-\frac{\eta}{2}+\epsilon}
\end{align*}
where we used the fact that $m\pm n$ have at most $(mn)^\epsilon$ divisors in $\Qcal$. If we use the bound \eqref{RS} after applying Cauchy--Schwarz inequality to the $m$-sum we get that the off-diagonal contribution is
	$O(Q^{-\frac{\eta}{2}+\epsilon})$.

The sums with the constraints $m\equiv \pm n[q_i]$ are bounded in the exact same way by $O(Q_i^{-1+\epsilon})$ (because of the factor $\phi(q_i)^{-1}$ in front of the diagonal terms).
To end the computation of $S_0$, we are left with the contribution of the trivial character. But once again, we have
\begin{align*}
	\frac{2}{|\Qcal|}\sum_{q\in \Qcal}\frac{1}{\phi(q)}& \sum_{\substack{m,n\geq1\\(mn,q_1q_2)=1 }}\frac{A(m,1)}{\sqrt{mn}}V\left(\frac{mn}{q^{2-\eta}}\right)\\
	&\ll Q^{-2+\epsilon}\sum_{q\in \Qcal} \sum_{mn < Q^{{2-\eta}+\epsilon}}|A(m,1)|(mn)^{-\frac{1}{2}+\epsilon} +O(Q^{-100})\\
	&\ll Q^{-\frac{\eta}{2}+\epsilon}.
\end{align*}
We used once again the bounds \eqref{boundV} for $V_i$, the Cauchy--Scwharz inequality and the Rankin--Selberg bound \eqref{RS}.

By putting all these estimates together, we get
\begin{align}
	S_0 = L(1,\pi) + O(Q^{-\frac{\eta}{2}+\epsilon}+Q^{\frac{\eta-2}{7}}+Q^{-\frac{9\delta}{14}+\epsilon}).\label{boundS0}
\end{align}

\subsection{Bound for the remaining sum $S_1$}

Le us rewrite the expression of $S_1$: 

\begin{align*}
	S_1 = \frac{2}{|\Qcal|}\sum_{q\in \Qcal}\frac{1}{\phi(q)}\sideset{}{^*}\sum_{\chi [q]}\epsilon(\chi)^2\sum_{m,n\geq1}\frac{A(1,m)\ol\chi(m)\chi(n)}{\sqrt{mn}}V\left(\frac{mn}{q^{2+\eta}}\right).
\end{align*}
In order to compute the contribution of $S_1$, we will need to cut the sum into dyadic intervals. Let $W$ be a smooth, compactly supported function on $[1,2]$. And define 
\begin{align*}
	\wt{S}(M,N) := \frac{2}{|\Qcal|}\sum_{q\in \Qcal}&\frac{1}{\phi(q)}\sideset{}{^*}\sum_{\chi [q]}\epsilon(\chi)^2\\&\times\sum_{m,n\geq1}\frac{A(1,m)\ol\chi(m)\chi(n)}{\sqrt{mn}}V\left(\frac{mn}{q^{2+\eta}}\right)W\left(\frac{m}{M}\right)W\left(\frac{n}{N}\right).
\end{align*}

Using the estimates on $V$ \eqref{boundV}, we can assume that $MN<Q^{2+\eta+\epsilon}$ for a cost of $O(Q^{-100})$. This implied that we only have to bound $O(\ln(Q))$ such sums. By expanding $V$ into an integral, we have
\begin{align*}
	\wt{S}(M,N) = \frac{1}{2\pi i}\int_{(\epsilon)}S(M,N)\frac{G\left(s+\frac{1}{2}\right)}{G\left(\frac{1}{2}\right)}\haar{s}
\end{align*}
where
\begin{align*}
	S(M,N) :=\frac{2}{|\Qcal|}\sum_{q\in \Qcal}\frac{q^{(2+\eta)s}}{\phi(q)}&\sideset{}{^*}\sum_{\chi [q]}\epsilon(\chi)^2\\&\times\sum_{m,n\geq1}\frac{A(1,m)\ol\chi(m)\chi(n)}{(mn)^{\frac{1}{2}+s}} W\left(\frac{m}{M}\right)W\left(\frac{n}{N}\right).
\end{align*}
Note that $S(M,N)$ only depends slightly on $s$ since $\Re{s} = \epsilon$.
We will treat two cases, depending on the relative size of $M$ and $N$.

\paragraph{Case 1: $N\geq MQ^{-1+\eta}$.} In this case, we will use Lemma \ref{Luo} after using a summation formula for the $n$-sum.
Using the Mellin inversion formula, we have\begin{align*}
	\sum_{n}\chi(n)n^{-\frac{1}{2}+s}W\left(\frac{n}{N}\right) = \frac{1}{2\pi i}\int_{(2)}N^{\omega}L\left(\frac{1}{2}+s+\omega,\chi\right)\wt W(\omega)d\omega
\end{align*}
where $\wt W$ is the Mellin transform of $W$, which is rapidly decreasing since $W$ is compactly supported. After shifting the line of integration to $\Re{\omega}=-2$, applying the functional equation of $L(s,\chi)$, doing the change of variable $
\omega\mapsto-\omega$ and expanding the L function,  we get 
\begin{align*}
	\sum_{n}\chi(n)n^{-\frac{1}{2}+s}W\left(\frac{n}{N}\right) = \epsilon(\chi)q^{-s}\sum_{n}\ol\chi(n)n^{-\frac{1}{2}-s}\Psi_1\left(\frac{nN}{q}\right)
\end{align*}
where
\begin{align*}
	\Psi_1(y) = \frac{1}{2\pi i}\int_{(2)}y^{-\omega}\frac{\Gamma\left(\frac{\omega-s+\frac{1}{2}}{2}\right)}{\Gamma\left(\frac{-\omega+s+\frac{1}{2}}{2}\right)}\wt{W}(-s)ds.
\end{align*}
By shifting the line of integration in the last expression to $\Re{s}=0$ and $\Re{s}= A$ for $A> 0$, we get
\begin{align*}
	\Psi_1(y) \ll \min(1,y^{-A}).
\end{align*}
In particular, we have 
\begin{align*}
	\sum_{n}\chi(n)n^{-\frac{1}{2}+s}W\left(\frac{n}{N}\right) =& \epsilon(\chi)q^{-s}\sum_{n<\frac{Q^{1+\epsilon}}{N}}\ol\chi(n)n^{-\frac{1}{2}-s}\Psi_1\left(\frac{nN}{q}\right) + O\left(Q^{-100}\right)\\
	=& \frac{1}{2\pi i}\int_{(0)}N^{-\omega}\epsilon(\chi)q^{-s+\omega}\\&\times\sum_{n<\frac{Q^{1+\epsilon}}{N}}\ol\chi(n)n^{-\frac{1}{2}-s+\omega}\frac{\Gamma\left(\frac{\omega-s+\frac{1}{2}}{2}\right)}{\Gamma\left(\frac{-\omega-s+\frac{1}{2}}{2}\right)}\wt{W}(-s)ds\\&+ O\left(Q^{-100}\right).
\end{align*}
By inserting this expression in the definition of $S(M,N)$, we only need to estimate the sums
\begin{align*}
	\sum_{k\geq1}\frac{a(k)}{k^\frac{1}{2}}B_3((1+\eta)s+\omega,k).
\end{align*}
with \begin{align*}
	B_3(s,k) = \frac{1}{|\mathcal{Q}|}\sum_{q\in\Qcal}\frac{q^{s}}{\phi(q)}\sideset{}{^*}\sum_{\substack{\chi [q_1]}}\epsilon(\ol{\chi})^3\chi(k)
\end{align*}
and \begin{align*}
	a(k) = \sum_{\substack{mn=k\\1\leq n< \frac{Q^{1+\epsilon}}{N}}}A(1,m)m^{-s}n^{-\omega+s}W\left(\frac{m}{M}\right).
\end{align*}
Here, we have $\Re{s}=\epsilon$ and $\Re{\omega}=0$. Using Cauchy's inequality, the above sum is bounded by 
\begin{align*}
	\left(\sum_{k\in \Z}\frac{|a(k)|^2}{k}H^{-1}\left(\frac{kN}{QM}\right)\right)^\frac{1}{2}\left(\sum_{k\in\Z}H\left(\frac{kN}{QM}\right)|B_3((1+\eta)s+\omega,k)|^2\right)^\frac{1}{2}
\end{align*}
where $H(y) = \pi^{-1}(1+y^2)^{-1}$.
Let us compute the first sum directly:
\begin{align*}
	\sum_{k\in \Z}\frac{|a(k)|^2}{k}H^{-1}\left(\frac{kN}{Qn}\right) &= Q^\epsilon\sum_{k\in \Z}k^{-1}\left|\sum_{\substack{mn=k\\1\leq n< \frac{Q^{1+\epsilon}}{N}\\m\asymp M}}|A(1,m)|\right|^2H^{-1}\left(\frac{kN}{QM}\right)\\
	&\ll Q^\epsilon\sum_{k\in \Z}k^{-1+\epsilon}\sum_{\substack{mn=k\\1\leq n< \frac{Q^{1+\epsilon}}{N}\\m\asymp M}} |A(1,m)|^2H^{-1}\left(\frac{kN}{QM}\right)\\
	&= Q^\epsilon\sum_{\substack{mn=k\\1\leq n< \frac{Q^{1+\epsilon}}{N}\\m\asymp M}} (mn)^{-1+\epsilon} |A(1,m)|^2H^{-1}\left(\frac{mnN}{QM}\right)\\
	&\ll Q^{\epsilon}.
\end{align*}
In order to estimate the square of the inner sum, we used the following bound
\begin{align*}
	\left|\sum_{n=1}^{N}c_n\right|^2\ll (2N-1)\sum_{n=1}^{N}|c_n|^2
\end{align*}
and noticed the fact that the sum had at most $d(k)\ll k^\epsilon$ terms. For the last inequality, we used once again the bound \eqref{RS} after doing a summation by parts.
The second sum is bounded by using Lemma \ref{Luo}:

\begin{align*}
	\sum_{k\in\Z}H\left(\frac{kN}{QM}\right)|B_3((1+\eta)s+\omega,k)|^2&\ll Q^{\epsilon}\left(\frac{QM}{N} Q^{-2+2\delta+\epsilon} +Q^{-\delta+\epsilon}\right)\\&\ll
	\frac{M}{N}Q^{-1+2\delta+\epsilon}+Q^{-\delta+ \epsilon}\\
	&\ll Q^{-\eta+2\delta+\epsilon}+Q^{-\delta+ \epsilon}.
\end{align*}
In the last line, we used the fact that $N\geq M Q^{-1+\eta}$. Also notice that $\delta$ is such that $\delta<1/4-3a/2$.
Finally, we have
\begin{align*}
	S(M,N)\ll Q^{-\frac{\eta}{2}+\delta+\epsilon}+Q^{-\frac{\delta}{2}+ \epsilon}.
\end{align*}

\paragraph{Case 2: $N< M Q^{-1+\eta}$.} In this case we will use a summation formula for the $m$-sum and we will use the hypothesis \eqref{hyp}. Note that this formula is similar to the usual Voronoï summation formula (see \cite{goldfeld2008voronoi} or \cite{miller2006automorphic}) (which is obtained after multiplying by $\epsilon(\chi)$ and summing over the charachers).
Proceeding as before, we have
\begin{align*}
	\sum_{m} A(1,m)\ol\chi(m)n^{-\frac{1}{2}-s}W\left(\frac{m}{M}\right) = \epsilon(\ol \chi)^3q^{-3s}\sum_{m}\frac{A(m,1)\chi(m)}{m^{\frac{1}{2}-s}}\Psi_2\left(\frac{mM}{q^3}\right)
\end{align*}
where
\begin{align*}
	\Psi_2(y) = \frac{1}{2\pi i}\int_{(2)}y^{-\omega}\frac{G_\pi\left(\omega -s+\frac{1}{2}\right)}{G_\pi\left(-\omega+s+\frac{1}{2}\right)}\wt{W}(-\omega)\haar{\omega}.
\end{align*}
We have once again 
\begin{align*}
	\Psi_2(y)\ll \min(1,y^{-A}).
\end{align*}
Inserting this into $S(M,N)$ gives
\begin{align*}
	S(M,N) = \frac{2}{|\Qcal|}\sum_{q\in \Qcal}\frac{q^{(-1+\eta)s}}{\phi(q)}&\sideset{}{^*}\sum_{\chi [q]}\epsilon(\ol\chi)\\&\times\sum_{m,n\geq1}\frac{A(m,1)\chi(mn)}{m^{\frac{1}{2}-s}n^{\frac{1}{2}+s}}\Psi_2\left(\frac{mM}{q^3}\right)W\left(\frac{n}{N}\right).
\end{align*}
We compute the character sum using the second part of Proposition \ref{characterSums}:

\begin{align*}
	S(M,N) = \frac{2}{|\Qcal|}\sum_{q\in \Qcal}&{q^{-\frac{1}{2}+(-1+\eta)s}}\sum_{(mn,q)=1}\frac{A(m,1)}{m^{\frac{1}{2}-s}n^{\frac{1}{2}+s}}\Psi_2\left(\frac{mM}{q^3}\right)W\left(\frac{n}{N}\right)\\ &\times\sum_{\alpha\in\{\pm1\}}\sum_{d|q}\frac{1}{\phi(d)}e_{\frac{q}{d}}(\alpha mn).
\end{align*}
In order to use the hypothesis \eqref{hyp}, we need to remove the condition $(mn,q)=1$. If $p$ is a prime number, define the function $\psi_p$ by 
\begin{align*}
	\psi_p(n) = 1 - \frac{1}{p}\sum_{b=0}^{p-1}e_p(bn)= \left\{\begin{aligned}
		 &1 \text{ if }(n,p)=1,\\
		 &0 \text{ otherwise}.
	\end{aligned}\right.
\end{align*}
By multiplying the expression we thus have
\begin{align*}
	S(M,N) = \frac{2}{|\Qcal|}\sum_{q\in \Qcal}&{q^{-\frac{1}{2}+(-1+\eta)s}}\sum_{n}n^{-\frac{1}{2}-s}W\left(\frac{n}{N}\right)\\&\times\sum_{m}\frac{A(m,1)}{m^{\frac{1}{2}-s}}\psi_2\left(\frac{mM}{q^3}\right) \\&\times\sum_{\alpha\in\{\pm1\}}\sum_{d|q}\frac{1}{\phi(d)}e_{\frac{q}{d}}(\alpha mn)\psi_{q_1}(mn)\psi_{q_2}(mn).
\end{align*}
The inner $m$-sum is a (normalized) linear combination of sums of the shape 
\begin{align*}
	\sum_m \frac{A(m,1)}{m^{\frac{1}{2}-s}}e(x m)\Psi_2\left(\frac{mM}{q^3}\right)
\end{align*}
for some real numbers $x$. Using the bound for $\Psi_2$, a summation by parts and the hypothesis \eqref{hyp}, these sums are uniformly bounded by $O\left(Q^{3a+\epsilon}M^{-a}\right)$ for some $a<1/6$.
This gives 
\begin{align*}
	S(M,N)\ll \frac{N^\frac{1}{2}}{M^a}Q^{3a-\frac{1}{2}+\epsilon}.
\end{align*}
Now, if $M\geq Q^{\frac{3}{2}}$, we use the bound $MN<Q^{2+\eta+\epsilon}$ to get
\begin{align*}
	S(M,N)\ll Q^{3a+\frac{1}{2}+\frac{\eta}{2}+\epsilon}M^{-\frac{1}{2}-a}\ll Q^{-\frac{1}{4}+\frac{3a}{2}+\frac{\eta}{2}+\epsilon}.
\end{align*}
If $M< Q^{\frac{3}{2}}$, we use the bound $N< M Q^{-1+\eta}$ to get once again
\begin{align*}
	S(M,N)\ll Q^{3a-1+\frac{\eta}{2}+\epsilon}M^{\frac{1}{2}-a}\ll Q^{-\frac{1}{4}+\frac{3a}{2}+\frac{\eta}{2}+\epsilon}.
\end{align*}

In conclusion, we have
\begin{align}\label{boundS1}
	S_1\ll Q^{-\frac{1}{4}+\frac{3a}{2}+\frac{\eta}{2}+\epsilon} + Q^{-\frac{\eta}{2}+\epsilon}+Q^{-\delta+ \epsilon}
\end{align}
\subsection{Conclusion of the proof}

Using the estimates \eqref{boundS0} and \eqref{boundS1},
we have 
\begin{align*}
	\Scal_\Qcal = L(1,\pi) + O\left(Q^{-\frac{\eta}{2}+\delta+\epsilon}+Q^{-\frac{1}{4}+\frac{3a}{2}+\frac{\eta}{2}+\epsilon}+Q^{\frac{\eta-2}{7}+\epsilon}+Q^{-\frac{\delta}{2}+ \epsilon}\right)
\end{align*}
Take for instance
\begin{align*}
	\eta = \frac{1}{4} - \frac{3a}{2}+\delta.
\end{align*}
The above estimate becomes
\begin{align*}
	\Scal_\Qcal = L(1,\pi) + O\left(Q^{-\frac{1}{8}+\frac{3a}{4}+\frac{\delta}{2}+\epsilon}+Q^{-\frac{\delta}{2}+ \epsilon}\right).
\end{align*}
Note that the best choice of $\delta$ is \begin{align*}
	\delta = \frac{1}{8}-\frac{3a}{4}
\end{align*}
which gives the estimate
\begin{align*}
		\Scal_\Qcal = L(1,\pi) + O\left(Q^{-\frac{1}{16}+\frac{3a}{8}+\epsilon}\right).
\end{align*}
However it is not in our interest to restrict the choice of $\delta$ since we are interested in the non vanishing of the moment for as many sets $\Qcal$ as possible.

We deduce immediately Corollary \ref{cor} because $\Scal_\Qcal$ is non zero for infinitely many disjoint sets $\Q_cal$, there exists infinitely many non zero terms in these sums.
Moreover,
\begin{align*}
	L\left(\frac{1}{2},\chi\right) = \ol{L\left(\frac{1}{2},\ol\chi\right)}.
\end{align*}

\bibliography{gl_non_vanishing}
\bibliographystyle{plain}

\end{document}